\documentclass[referee]{gji}
\usepackage{timet}
\usepackage{amsmath, amssymb}
\usepackage{enumerate}
\usepackage{graphicx}
\graphicspath{{fig/}}

\newtheorem{example}{Example}[section]
\newtheorem{remark}{Remark}[section]

\DeclareMathOperator*{\argmin}{argmin}

\newcommand{\br}{\mathbb{R}}
\newcommand{\mca}{\mathcal{A}}
\newcommand{\mcr}{\mathcal{R}}
\newcommand{\rd}{\mathrm{d}}
\newcommand{\veps}{\varepsilon}

\newcommand{\bn}{\boldsymbol{n}}
\newcommand{\boeta}{\boldsymbol{\eta}}
\newcommand{\bx}{\boldsymbol{x}}
\newcommand{\bxi}{\boldsymbol{\xi}}

\newcommand{\abs}[1]{\left\vert#1\right\vert}
\newcommand{\brac}[1]{\left(#1\right)}
\newcommand{\norm}[1]{\left\Vert#1\right\Vert}

\title[A new earthquake location method based on the waveform inversion]{A new earthquake location method based on the waveform inversion}
\author[H. Wu et al.]
	{Hao Wu$^1$, Jing Chen$^1$, Xueyuan Huang$^2$, and Dinghui Yang$^1$\\
	$^1$ Department of Mathematical Sciences, Tsinghua University, Beijing, China. E-mail: dhyang@math.tsinghua.edu.cn \\
	$^2$ Department of Mathematics, School of Science, Beijing Technology and Business University, Beijing, China.}
\date{}
\pagerange{\pageref{firstpage}--\pageref{lastpage}}
\volume{}
\pubyear{}

\begin{document}

\label{firstpage}

\maketitle

\begin{summary}
	In this paper, a new earthquake location method based on the waveform inversion is proposed. As is known to all, the waveform misfit function is very sensitive to the phase shift between the synthetic waveform signal and the real waveform signal. Thus, the convergence domain of the conventional waveform based earthquake location methods is very small. In present study, by introducing and solving a simple sub-optimization problem, we greatly expand the convergence domain of the waveform based earthquake location method. According to a large number of numerical experiments, the new method expands the range of convergence by several tens of times. This allows us to locate the earthquake accurately even from some relatively bad initial values.
\end{summary}



\section{Introduction} \label{sec:intro}

The earthquake location is a fundamental problem in seismology \cite{Ge:03,Th:14}. It consists of two parts: the determination of hypocenter $\bxi$ and origin time $\tau$. These information are extremely important in quantitative seismology, e.g. the earthquake early warning system \cite{SaLoZo:08}, the investigation of seismic heterogeneous structure \cite{ToYaLiYaHa:16,WaEl:00}. In particular, there are also significant interests in micro-earthquake which has many applications in exploration seismology \cite{LeSt:81,PrGe:88}.

Due to the importance of the earthquake location problem, numerous studies have been done theoretically and experimentally, see e.g. \cite{Ge:03,Ge:03b,Ge:10,Ge:12,Mi:86,PrGe:88,Th:85}. However, many studies are based on the ray theory, which has low accuracy when the wave length is not small enough compared to the scale of wave propagating region \cite{EnRu:03,JiWuYa:08,JiWuYa:11,LiRuTa:13,RaPoFi:10,WuYa:13}.  This may lead to inaccurate or even incorrect earthquake location results. An alternative way is to solve the wave equation directly to get accurate information for inversion. This method is becoming popular in recent years, as a result of the fast increasing computational power \cite{HuYaToBaQi:16,KiLiTr:11,LiGu:12,LiPoKoTr:04,TaLiTr:07,ToChKoBaLi:14,ToZhYaYaChLi:14,ToZhYaYaChLi:14b}.

In the work by \cite{LiPoKoTr:04}, see also \cite{KiLiTr:11}, the spectral-element solvers are implemented to invert the basic information of earthquakes. The misfit functions defined based upon the envelope of the waveforms are minimized to provide the best estimation of source model parameters. Another approach proposed by \cite{ToYaLiYaHa:16}  is based on the wave-field relation between the hypocenter $\bxi$ and its perturbation $\bxi+\delta \bxi$  \cite{Al:10}. Due to the foregoing observation, the travel-time differences between the synthetic single and the real single can be approximately expressed as the linear function of hypocenter perturbation $\delta \bxi$. The authors then derived the sensitivity kernel by using the forward and adjoint wavefields.

However, the above mentioned papers on the earthquake location are not directly used the waveform difference since the waveform misfit function is very sensitive to the phase shift between the synthetic waveform signal and the real waveform signal \cite{LiPoKoTr:04}. Consider the bad mathematical properties of the delta function $f(t-\tau)\delta(\bx-\bxi)$, who is  appeared as the source of wave equation, even small perturbation of hypocenter $\delta \bxi$ and origin time $\delta \tau$ would generate large deviation of waveform. Thus, it is not surprising that the range of convergence of the conventional waveform based method is very small. On the other hand, the waveform signal may contain more information, which could lead to more accurate location result. Thus, it is necessary to develop new techniques to expand the convergence domain of the waveform based location method. 

In this paper, we present a new method to locate the earthquake accurately. For the sake of simplicity, we use the acoustic wave equation and only deal with the earthquake hypocenter and origin time. There is no essential difficulty to consider the elastic wave equation or involve more earthquake information, e.g. the moment magnitudes \cite{LiPoKoTr:04}. The starting point is to keep $\frac{\norm{\delta s(\bx,t)}}{\norm{s(\bx,t)}}\ll1$ in a modified sense. This is a fundamental assumption of the first-order Born approximation in the adjoint method. But it is not easy to guarantee in the classical sense, even if $\frac{\norm{\delta \bxi}}{\norm{\bxi}}$ and $\frac{\norm{\delta \tau}}{\norm{\tau}}$ are small. This is due to the bad mathematical properties of the delta function $f(t-\tau)\delta(\bx-\bxi)$ in the wave equation. To solve this problem, we shift the synthetic data so that its difference with the real data is minimized. The shifting parameter can be obtained by solving a simple sub-optimization problem. The above effects ensure correctness of the important assumption $\frac{\norm{\delta s(\bx,t)}}{\norm{s(\bx,t)}}\ll1$ of the adjoint method in a large range. Thus, we can expect a large convergence domain of the new earthquake location method. According to the numerical experiments, the range of convergence is significantly enlarged. We also remark that the idea is similar to the Wasserstein metric \cite{EnFr:14,EnFrYa:16}, but we provide a simple and alternative implementation.

The paper is organized as follows. In Section \ref{sec:jim}, the conventional waveform based adjoint inversion method is reviewed for the earthquake hypocenter and origin time. We propose the new method for the earthquake location in Section \ref{sec:iie}. In Section \ref{sec:num}, the numerical experiments are provided to demonstrate the effectiveness of the new method. Finally, we make some conclusive remarks in Section \ref{sec:con}.

\section{The inversion method} \label{sec:jim}
Consider the scalar acoustic wave equation
\begin{equation} \label{eqn:wave}
	\frac{\partial^2 u(\bx,t;\bxi,\tau)}{\partial t^2}=\nabla \cdot\brac{c^2(\bx)\nabla u(\bx,t;\bxi,\tau)}
		+f(t-\tau)\delta(\bx-\bxi), \quad \bx,\bxi\in\Omega,
\end{equation}
with initial-boundary conditions
\begin{align}
	& u(\bx,0;\bxi,\tau)=\partial_t u(\bx,0;\bxi,\tau)=0, \quad \bx\in\Omega, \label{ic:wave} \\
	& \bn\cdot\brac{c^2(\bx)\nabla u(\bx,t;\bxi,\tau)}=0, \quad \bx\in\partial\Omega. \label{bc:wave}
\end{align}
Here $u(\bx,t;\bxi,\tau)$ is the wavefield with respect to parameters $\tau$ and $\bxi$. The wave speed is $c(\bx)$. The simulated domain $\Omega\subset \br^d$, $d$ is the dimension of the problem and $\bn$ is the unit outer normal vector to the boundary $\partial \Omega$ of $\Omega$. The seismogram at source has the form of Ricker wavelet
\begin{equation*}
	f(t)=A\brac{1-2\pi^2f_0^2t^2}e^{-\pi^2f_0^2t^2},
\end{equation*}
in which $f_0$ is the dominant frequency and $A$ is the normalization factor. In this study, the point source hypothesis $\delta(\bx-\bxi)$ for the hypocenter focus is considered for the situation where the temporal and spatial scales of seismic rupture are extremely small compared to the scales of seismic waves propagated \cite{AkRi:80,Ma:15}. For simplicity, the reflection boundary condition \eqref{bc:wave} is considered here. There is no essential difference for other boundary conditions, e.g. the perfectly matched layer absorbing boundary condition \cite{KoTr:03,MaYaSo:15}.

\begin{remark} \label{rem:time_tran}
	For the acoustic wave equation \eqref{eqn:wave}, we have the invariance property in time translation \cite{Ev:10}
	\begin{equation*}
		u(\bx,t-\Delta \tau;\bxi,\tau)=u(\bx,t;\bxi,\tau+\Delta \tau). \; \Box
	\end{equation*}
\end{remark}

\begin{remark}
	The compatibility condition requires that
	\begin{equation*}
		f(t-\tau)=0 \;\textrm{and}\; \tau>0.
	\end{equation*}
	These are very nature in practical problems. $\Box$
\end{remark}

Let $\bxi_T$ and $\tau_T$ be the real earthquake hypocenter and origin time. Thus, the real earthquake signal $d_r(t)$, which was receiver at station $r$ can be considered as
\begin{equation} \label{eqn:d_rel}
	d_r(t)=u(\boeta_r,t;\bxi_T,\tau_T).
\end{equation}
Here $\boeta_r$ is the location of the $r-$th receiver. The synthetic signal $s(\bx,t)$ corresponding to the initial hypocenter $\bxi$ and origin time $\tau$ is
\begin{equation} \label{eqn:s_rel}
	s(\bx,t)=u(\bx,t;\bxi,\tau).
\end{equation}
By introducing the misfit function
\begin{equation} \label{fun:chi}
	 \chi_r(\bxi,\tau)=\frac{\int_0^T\abs{d_r(t)-s(\boeta_r,t)}^2\rd t}{2\int_0^T\abs{d_r(t)}^2\rd t},
\end{equation}
we define the nonlinear optimization problem
\begin{equation} \label{prob:non_opti}
	(\bxi_T,\tau_T)= \mathop{\argmin}_{\bxi,\;\tau}\sum_r\chi_r(\bxi,\tau).
\end{equation}
Obviously, the global solution exists and is unique \cite{NoWr:99}. In the following part, the sensitivity kernel \cite{LiGu:12,RaPoFi:10,ToZhYaYaChLi:14} will be derived to solve this inversion problem iteratively.

\subsection{The adjoint method} \label{subsec:adj_method}
The perturbation of parameters $\frac{\norm{\delta \bxi}}{\norm{\bxi}}$ and $\frac{\norm{\delta \tau}}{\norm{\tau}}\ll1$ would generate the perturbation of wave function $\delta s(\bx,t)$, it writes
\begin{equation*}
	\delta s(\bx,t)=u(\bx,t;\bxi+\delta \bxi,\tau+\delta \tau)-u(\bx,t;\bxi,\tau).
\end{equation*}
Then $\delta s(\bx,t)$ satisfies
\begin{equation}  \label{eqn:adj_sub}
	\left\{\begin{array}{ll}
		\frac{\partial^2 \delta s(\bx,t)}{\partial t^2}=\nabla\cdot\brac{c^2(\bx)\nabla \delta s(\bx,t)}
			+f(t-(\tau+\delta \tau))\delta(\bx-(\bxi+\delta \bxi))-f(t-\tau)\delta(\bx-\bxi), & \bx\in\Omega, \\
                \delta s(\bx,0)=\frac{\partial \delta s(\bx,0)}{\partial t}=0, & \bx\in\Omega, \\
                \bn\cdot\brac{c^2(\bx)\nabla \delta s(\bx,t)}=0, & \bx\in\partial\Omega.
        \end{array}\right.
\end{equation}
Multiply an arbitrary test funciton $w_r(\bx,t)$ on equation \eqref{eqn:adj_sub}, integrate it on $\Omega\times[0,T]$  and use integration by parts, we obtain
\begin{multline} \label{eqn:integration}
	\int_0^T\int_{\Omega}\frac{\partial^2 w_r}{\partial t^2}\delta s\rd \bx\rd t
		-\int_{\Omega}\left.\frac{\partial w_r}{\partial t}\delta s\right|_{t=T}\rd \bx
		+\int_{\Omega}\left.w_r\frac{\partial \delta s}{\partial t}\right|_{t=T}\rd \bx \\
	=\int_0^T\int_{\Omega}\delta s\nabla\cdot(c^2\nabla w_r)\rd \bx \rd t
		-\int_0^T\int_{\partial \Omega}\bn\cdot (c^2\nabla w_r)\delta s\rd \zeta \rd t \\
                 +\int_0^T f(t-(\tau+\delta \tau))w_r(\bxi+\delta \bxi,t)-f(t-\tau)w_r(\bxi,t)\rd t \\
	\approx\int_0^T\int_{\Omega}\delta s\nabla\cdot(c^2\nabla w_r)\rd \bx \rd t
                  -\int_0^T\int_{\partial \Omega}\bn\cdot (c^2\nabla w_r)\delta s\rd \zeta \rd t \\
                  +\int_0^T f(t-\tau)\nabla w_r(\bxi,t)\cdot \delta \bxi-f'(t-\tau)w_r(\bxi,t)\delta \tau\rd t.
\end{multline}
Note that the Taylor expansion is used and higher order terms is ignored in the last step.

On the other hand, the misfit function \eqref{fun:chi} also generates the perturbation with respect to $\delta s(\bx,t)$, assume that $\frac{\norm{\delta s(\bx,t)}}{\norm{s(\bx,t)}}\ll 1$, it writes
\begin{multline} \label{fun:deltachi}
	\delta \chi_r=\chi_r(\bxi+\delta \bxi,\tau+\delta \tau)-\chi_r(\bxi,\tau) \\
	=\frac{\int_0^T\brac{\abs{d_r(t)-(s+\delta s)(\boeta_r,t)}^2-\abs{d_r(t)-s(\boeta_r,t)}^2}\rd t}{2\int_0^T\abs{d_r(t)}^2\rd t} \\
		\approx-\frac{\int_0^T(d_r(t)-s(\boeta_r,t))\delta s(\boeta_r,t)\rd t}{\int_0^T\abs{d_r(t)}^2\rd t} \\
	=-\frac{\int_0^T\int_{\Omega}(d_r(t)-s(\boeta_r,t))\delta s(\bx,t)\delta(\bx-\boeta_r)\rd \bx\rd t}{\int_0^T\abs{d_r(t)}^2\rd t},
\end{multline}
where ``$\approx$'' is obtained by ignoring high order terms of $\delta s(\bx,t)$.

Let $w_r(\bx,t)$ satisfy the wave equation with terminal-boundary conditions
\begin{equation}
	\left\{\begin{array}{ll}
		\frac{\partial^2 w_r(\bx,t)}{\partial t^2}=\nabla\cdot\brac{c^2(\bx)\nabla w_r(\bx,t)}
			+\frac{d_r(t)-s(\boeta_r,t)}{\int_0^T\abs{d_r(t)}^2\rd t}\delta(\bx-\boeta_r), & \bx\in\Omega, \\
                 w_r(\bx,T)=\frac{\partial w_r(\bx,T)}{\partial t}=0, & \bx\in\Omega, \\
                 	 \bn\cdot\brac{c^2(\bx)\nabla w_r(\bx,t)}=0, & \bx\in\partial\Omega.
	\end{array}\right.
\end{equation}
Thus, the linear relation for $\delta \chi$ and $\delta \bxi, \;\delta \tau$ can be obtained by subtracting \eqref{fun:deltachi} from \eqref{eqn:integration} 
\begin{equation} \label{eqn:delchi_rel}
	-\delta \chi_r=\int_0^T f(t-\tau)\nabla w_r(\bxi,t)\cdot \delta \bxi-f'(t-\tau)w_r(\bxi,t)\delta \tau\rd t.
\end{equation}
In particular, if
\begin{equation*}
	\bxi+\delta \bxi=\bxi_T,\; \tau+\delta \tau=\tau_T,
\end{equation*}
it implies
\begin{equation*}
	\chi_r(\bxi+\delta \bxi,\tau+\delta \tau)=0\; \Rightarrow\; \delta \chi_r=-\chi_r(\bxi,\tau).
\end{equation*}
This gives an alternative form of equation \eqref{eqn:delchi_rel}
\begin{equation} \label{eqn:chi_rel}
	\chi_r(\bxi,\tau)=\int_0^T f(t-\tau)\nabla w_r(\bxi,t)\cdot \delta \bxi-f'(t-\tau)w_r(\bxi,t)\delta \tau\rd t.
\end{equation}

By defining the sensitivity kernel for the hypocenter $\bxi$ and origin time $\tau$ as
\begin{align*}
         & K^{\bxi}_r = \int_0^T \nabla w_r(\bxi,t)f(t-\tau)\rd t, \\
         & K^{\tau}_r = -\int_0^T w_r(\bxi,t)f'(t-\tau)\rd t,
\end{align*}
equation \eqref{eqn:chi_rel} gives a single equation of the large linear system
\begin{equation} \label{eqn:large_linear_system}
	\frac{K^{\bxi}_r}{\chi_r(\bxi,\tau)}\cdot\delta \bxi+\frac{K^{\tau}_r}{\chi_r(\bxi,\tau)}\delta \tau=1.
\end{equation}
The above linear system has been nomailized so that the condition number can be optimized

\section{A new method to expand the convergence domain} \label{sec:iie}
In this section, we are investigating the techniques to enlarge the convergence domain for the inversion of earthquake hypocenter $\bxi_T$ and origin time $\tau_T$. It is assumed that the wave speed $c(\bx)$ is already well known. For situations of inaccurate or unknown wave speed, we refer to the discussions in \cite{LiPoKoTr:04} or the joint inversion for wave speed, hypocenter and origin time. 

\subsection{Estimation of the origin time} \label{subsec:est_ori_time}
As it was discussed in Section \ref{subsec:adj_method}, the first-order Born approximation in the adjoint method requires an infinitesimal perturbation assumption of wave function 
\begin{equation*}
	\frac{\norm{\delta s(\bx,t)}}{\norm{s(\bx,t)}}\ll 1,
\end{equation*}
see also \cite{LiGu:12,RaPoFi:10,ToZhYaYaChLi:14,TrTaLi:05}. However, as we will see in Example \ref{exam:large_pert}, it is very difficult to guarantee this assumption even if the perturbations of the earthquake hypocenter and origin time are very small
\begin{equation*}
	\frac{\norm{\delta \bxi}}{\norm{\bxi}}\ll1 \;\textrm{ and }\;  \frac{\norm{\delta \tau}}{\norm{\tau}}\ll1.
\end{equation*}
That's one of the reasons why the convergence domain of the  waveform based method is very small.

\begin{figure} 
	\centering
	\includegraphics[width=0.76\textwidth, height=0.30\textheight]{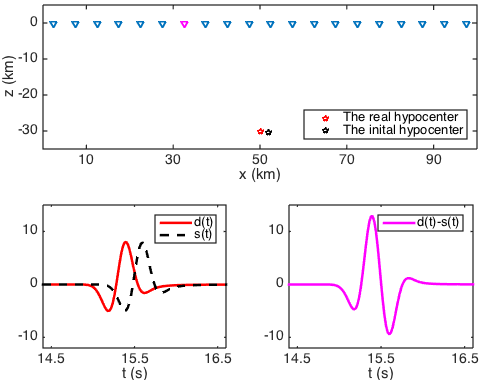}
	\caption{Illustration of Exam \ref{exam:large_pert}. Up: The receiver stations (up triangle, magenta for $r=7$ and blue for others), the real (red pentagram) and initial hypocenter (black pentagram). Bottom Left: The real signal $d_r(t)$ (red solid line) and the synthetic signal $s(\boeta_r,t)$ (black dashed line) at receiver station $r=7$. Bottom Right: The difference between the real and synthetic signal $d_r(t)-s(\boeta_r,t)$ (magenta solid line) at receiver station $r=7$. The text representation in the figure has been simplified without causing any misunderstandings.} \label{fig:large_pert}
\end{figure}

\begin{example} \label{exam:large_pert}
	This is a 2D unbounded problem with constant wave speed $c(\bx)\equiv c_0$ for the scalar acoustic wave equation \eqref{eqn:wave} with initial condition \eqref{ic:wave}. Its solution can be analytically given \cite{Ev:10}:
	\begin{equation} \label{sol:analy_2D}
		u(\bx,t;\bxi,\tau)=\left\{\begin{array}{ll}
			\frac{1}{2\pi c_0^2}\int_0^{\theta_0}\frac{f(\theta-\tau)}{\sqrt{(t-\theta)^2-(t-\theta_0)^2}}\rd \theta, & \theta_0>0, \\
			0, & \theta_0\le0, 
		\end{array}\right.
	\end{equation}
	in which
	\begin{equation*}
		\theta_0=t-\frac{1}{c_0}\norm{\bx-\bxi}_2.
	\end{equation*}
	Let $\bx=(x,z)$ denote the horizontal and depth coordinate respectively. The constant wave speed is $c_0=6.5km/s$. There are 20 equidistant receivers on the surface,
	\begin{equation*}
		\boeta_r=(x_r,z_r)=(5r-2.5km,0), \quad r=1,2,\cdots,20.
	\end{equation*}
	Consider an earthquake occurs at hypocenter $\bxi_T=(50km,-30km)$ and origin time $\tau_T=10s$ with dominant frequency $f_0=2$Hz, its signal $d_r(t)$ received by station $r=7$ can be considered as \eqref{eqn:d_rel} and \eqref{sol:analy_2D}. The synthetic signal $s(\boeta_r,t)$ corresponding to the initial hypocenter $\bxi=(52km,-30.3km)$ and origin time $\tau=10s$ at the same station $r=7$ can be obtained by \eqref{eqn:s_rel} and \eqref{sol:analy_2D}. The perturbation between the real and initial hypocenter is small
	\begin{equation*}
		\delta \bxi=\bxi_T-\bxi=(-2km,0.3km),
	\end{equation*}
	and it is also correct for the perturbation between the real and initial origin time
	\begin{equation*}
		\delta \tau=\tau_T-\tau=0.
	\end{equation*}
	Nevertheless, as we can see in Figure \ref{fig:large_pert}, the difference between the real signal $d_r(t)$ and the synthetic signal $s(\boeta_r,t)$ at receiver station $r=7$ is significantly large:
	\begin{equation*}
		\frac{\norm{d_r(t)-s(\boeta_r,t)}}{\norm{d_r(t)}}\sim 1,
	\end{equation*}
	which contracts to the basic assumption. $\Box$
\end{example}

The key observation in Exam \ref{exam:large_pert} is that the infinitesimal perturbation assumption $\frac{\norm{\delta s(\bx,t)}}{\norm{s(\bx,t)}}\ll 1$ is not trivial to get. However, we note that the main difference between the real signal $d_r(t)$ and synthetic signal $s(\boeta_r,t)$ is caused by the time shift. Thus, we define the relative error function with respect to the time shift of the synthetic signal
\begin{equation} \label{eqn:def_errfun}
	e_r(\tau)=\frac{\norm{d_r(t)-s(\boeta_r,t-\tau)}}{\norm{d_r(t)}}.
\end{equation}
Solving the following sub-optimization problem
\begin{equation} \label{eqn:comp_tau}
	\tau_r^*=\argmin_{\tau}e_r(\tau),
\end{equation}
the infinitesimal perturbation assumption may be satisfied in the sense of time translation.
\begin{equation} \label{eqn:new_assump}
	\frac{\norm{d_r(t)-s(\boeta_r,t-\tau_r^*)}}{\norm{d_r(t)}}\ll1,
\end{equation}

\begin{figure} 
	\centering
	\includegraphics[width=0.76\textwidth, height=0.30\textheight]{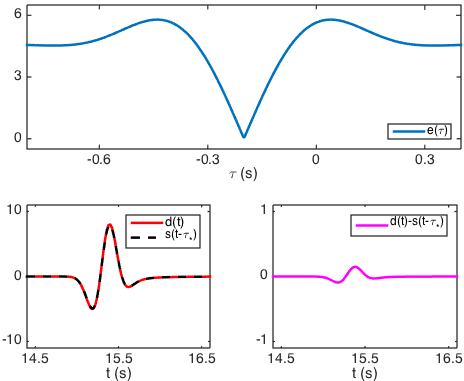}
	\caption{Illustration of Exam \ref{exam:shift_small_pert}. Up: The relative error function $e_r(\tau)$ defined in \eqref{eqn:def_errfun}. Bottom Left: The real signal $d_r(t)$ (red solid line) and the shifted synthetic signal $s(\boeta_r,t-\tau_r^*)$ (black dashed line) at receiver station $r=7$. Bottom Right: The difference between the real signal and the shifted synthetic signal $d_t(t)-s(\boeta_r,t-\tau_r^*)$ (magenta solid line) at receiver station $r=7$. The text representation in the figure has been simplified without causing any misunderstandings.} \label{fig:shift_small_pert}
\end{figure}

\begin{example} \label{exam:shift_small_pert}
	Consider the same parameters set up as in Exam \ref{exam:large_pert}, thereby the real signal $d_r(t)$ and the synthetic signal $s(\boeta_r,t)$ are the same as those in Exam \ref{exam:large_pert}.
	
	The relative error function $e_r(\tau)$ defined in \eqref{eqn:def_errfun} is presented in Figure \ref{fig:shift_small_pert} Up. We can observe a global minimum of $e_r(\tau)$. Thus, the optimal time translation parameter $\tau_r^*$ can be easily computed through \eqref{eqn:comp_tau}.
	
	According to the above time translation, the difference between the real signal $d_r(t)$ and the shifted synthetic signal $s(\boeta_r,t-\tau_r^*)$ is small, see Figure \ref{fig:shift_small_pert} Bottom. This implies that the modified infinitesimal perturbation assumption in \eqref{eqn:new_assump} is satisfied here. $\Box$
\end{example}

At last, by the invariance property in time translation, see Remark \ref{rem:time_tran}, this optimal time shift $\tau_r^*$ computed from \eqref{eqn:comp_tau} can be used to shift the initial origin time
\begin{equation} \label{eqn:new_initial_oritime}
	\hat{\tau}=\tau+\tau_r^*,
\end{equation}
so that the infinitesimal perturbation assumption $\frac{\norm{\delta s(\bx,t)}}{\norm{s(\bx,t)}}\ll 1$ can be satisfied in the original sense rather than the modified sense \eqref{eqn:new_assump}.

\begin{figure} 
	\centering
	\includegraphics[width=0.76\textwidth, height=0.30\textheight]{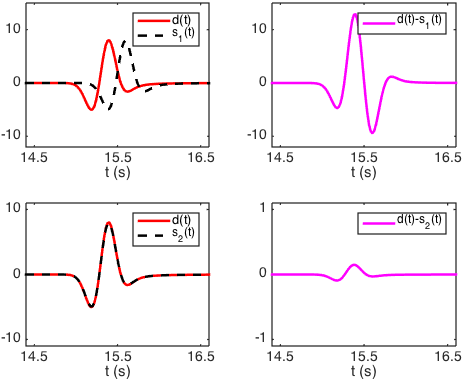}
	\caption{Illustration of Exam \ref{exam:comp_pert}. Left: The real signal $d_r(t)$ (red solid line) and the synthetic signal $s_i(t,\boeta_r)$ (black dashed line) at receiver station $r=7$. Right: The difference between the real signal and the synthetic signal $d_r(t)-s_i(\boeta_r,t)$ (magenta solid line) at receiver station $r=7$. Up: $i=1$, which corresponding to the initial origin time $\tau_1$. Bottom: $i=2$, which correspoding to the other initial origin time $\tau_2$. The text representation in the figure has been simplified without causing any misunderstandings.} \label{fig:comp_pert}
\end{figure}

\begin{example} \label{exam:comp_pert}
	Consider the same parameters set up as in Exam \ref{exam:large_pert}, thereby the real signal $d_r(t)$ is the same as in Exam \ref{exam:large_pert}. The synthetic signals are corresponding to the initial hypocenter $\bxi=(52km,-30.3km)$ and two different initial origin time: (1) $\tau_1=10$s, (2) $\tau_2=\tau_1+\tau_r^*$. We still focus on the signals received at station $r=7$. These synthetic signals can be obtained by
	\begin{equation*}
		s_1(\boeta_r,t)=u(\boeta_r,t;\bxi,\tau_1), \quad s_2(\boeta_r,t)=u(\boeta_r,t;\bxi,\tau_2).
	\end{equation*}
	In Figure \ref{fig:comp_pert} Up, the difference between the real signal $d_r(t)$ and the synthetic signal $s_1(\boeta_r,t)$ is large, but the difference between the real signal $d_r(t)$ and the other synthetic signal $s_2(\boeta_r,t)$ is small, see Figure \ref{fig:comp_pert} Bottom. $\Box$
\end{example}

\subsection{The selection of receiver stations} \label{subsec:rec}
In previous subsection, the time shift $\tau_r^*$ has been discussed for single receiver station $r$. For practical problems, there are many receiver stations, thus we need to solve the following sub-optimization problem
\begin{equation} \label{eqn:comp_gentau}
	\tau^*=\argmin_{\tau}\sum_{r\in\mcr}e_r(\tau),
\end{equation}
rather than \eqref{eqn:comp_tau}. Here $\mcr$ is the set of all receiver stations that we use for inversion, which will be determined later. The set of all receiver stations is denote by $\mca$, and it is obviously that $\mcr\subset\mca$.

\begin{figure} 
	\centering
	\includegraphics[width=0.76\textwidth, height=0.30\textheight]{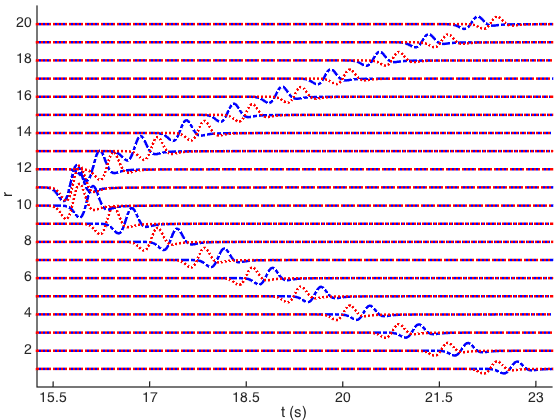}
	\caption{Illustration of Exam \ref{exam:diffrec_timeser}. The waveform comparison of real signals $d_r(t)$ (red dotted line) and the synthetic signals $s(\boeta_r,t)$ (blue dashdot line) for all receivers $r=1,2,\cdots,20$. The horizontal axis is the time $t$, and the longitudinal axis is the index of receiver $r$.} \label{fig:diffrec_timeser}
\end{figure}

\begin{example} \label{exam:diffrec_timeser}
	Consider the same parameters set up as in Exam \ref{exam:large_pert}, thereby the real signals $d_r(t)$ and the synthetic signals $s(\boeta_r,t)$ can be obtained in the same manner as in Exam \ref{exam:large_pert} 
	\begin{equation*}
		d_r(t)=u(\boeta_r,t;\bxi_T,\tau_T), \quad s(\boeta_r,t)=u(\boeta_r,t;\bxi,\tau)
	\end{equation*}
	for all receiver stations $r\in\mca=\{1,2,\cdots,20\}$. In Figure \ref{fig:diffrec_timeser}, we output all the real signals $d_t(t)$ and the synthetic signals $s(\boeta_r,t)$. According to the figures, we can see that
	\begin{align*}
		& \tau_r^*<0, \;\textrm{for}\; r=1,2,\cdots,10, \\
		& \tau_r^*>0, \;\textrm{for}\; r=11,12,\cdots,20.
	\end{align*}
	Therefore, we cannot get satisfactory value $\tau^*$ from \eqref{eqn:comp_gentau} if $\#\mcr$ is large. Here $\#\mcr$ denotes the number of elements in the set $\mcr$. $\Box$
\end{example}

The above example shows that discussions in Subsection \ref{subsec:est_ori_time} may fail when $\#\mcr$ is large. In fact, due to the small degree of freedom of the earthquake location problem, it is not necessary to consider large $\#\mcr$. On the other hand, $\#\mcr$ is the number of constraints, which is proportional to the number of wave field computations. Therefore, we prefer to choose a relative small $\#\mcr$ for inversion. Accordingly to our numerical experiences, a suitable choice of $\#\mcr$ is $5\sim7$. However, this discussion doesn't determine which elements should be in $\mcr$. A natural consideration is to solve a more general nonlinear optimization problem 
\begin{equation} \label{eqn:general_tau}
	(\tau^*,\mcr^*)=\argmin_{\tau,\;\mcr\subset\mca}\sum_{r\in\mcr}e_r(\tau).
\end{equation}
The essence of the above problem is that receivers set $\mcr$ is considered as optimization variable. It is easy to check that for $1\le n_1<n_2\le \#\mca$, we have
\begin{equation*}
	\min_{\tau,\;\#\mcr=n_1,\;\mcr\subset\mca}\sum_{r\in\mcr}e_r(\tau)
		\le \min_{\tau,\;\#\mcr=n_2,\;\mcr\subset\mca}\sum_{r\in\mcr}e_r(\tau).
\end{equation*}
In practice, solving the problem \eqref{eqn:general_tau} is complicated. Instead, we can firstly solve a simplified optimization problem
\begin{equation} \label{eqn:sim_optim}
	(\bar{\tau},\mcr^*)=\argmin_{\tau,\;\mcr\subset\mca}\sum_{r\in\mcr}\abs{\tau_r^*-\tau}^2.
\end{equation}
Then, for fixed receivers set $\mcr^*$, we have
\begin{equation} \label{eqn:final4tau}
	\tau^*=\argmin_{\tau}\sum_{r\in\mcr^*}e_r(\tau).
\end{equation}
Similar to equation \eqref{eqn:new_initial_oritime}, the optimal time shift $\tau^*$ for multiple receivers can also be used to shift the initial origin time
\begin{equation} \label{eqn:new_mult_initial_oritime}
	\hat{\tau}=\tau+\tau^*.
\end{equation}

\subsection{The detailed implementation} \label{subsec:imple}
In summary of all the above, the detailed implementation of the algorithm is as follows:
\begin{enumerate}[1.]
	\item Initialization. Set the tolerance value $\veps=0.01km$, the threshold value $\sigma=100km$ and the break-off step $K=30$. Let $k=0$ and give the initial hypocenter $\bxi_0$ and the initial origin time $\tau_0=0$. 
	
	\item For $\bxi_k$, solving \eqref{eqn:sim_optim} to determine the receivers set $\mcr_k^*$ and estimating the origin time $\tau_k$ by \eqref{eqn:final4tau} and \eqref{eqn:new_mult_initial_oritime}.
	
	\item Construct the sensitivity kernels $K_{r,k}^{\bxi},\;K_{r,k}^{\tau}$ for $r\in\mcr_k^*$ and solve the normalized linear system \eqref{eqn:large_linear_system} to get $\delta \bxi_k$ and $\delta \tau_k$, then update the estimation of hypocenter for step $k+1$,
		\begin{equation*}
			\bxi_{k+1}=\bxi_k+\delta \bxi_k.
		\end{equation*}
		
	\item If $\norm{\bxi_k-\bxi_{k+1}}<\veps$, go to step 7; If $\norm{\bxi_k-\bxi_{k+1}}>\sigma$, go to step 6.
	
	\item If $k+1>K$, go to step 6. Otherwise, let $k=k+1$ and go to step 2 for another iteration.
	
	\item Output the error message: ``The iteration diverges.'' and stop.
	
	\item	Update the estimation of origin time for step $k+1$,
		\begin{equation*}
			\tau_{k+1}=\tau_k+\delta \tau_k.
		\end{equation*}
		Output $(\bxi_{k+1},\tau_{k+1})$ and stop.$\;\Box$
\end{enumerate}

Once the value $(\bxi_{k+1},\tau_{k+1})$ is output, we get the hypocenter and the origin time for the specific earthquake. Otherwise, the algorithm should be restarted with different initial value of hypocenter $\bxi_0$ until the convergent result is obtained.

In this algorithm, the extra computational cost arise from solving the sub-optimization problem \eqref{eqn:sim_optim} and \eqref{eqn:final4tau}. But this part in the overall computational cost is minor. The reason is that the sub-optimization problem \eqref{eqn:sim_optim} and \eqref{eqn:final4tau} are only one dimensional. Taking into consideration the saving from less computation of the wave equations, the total cost is reduced here. Furthermore, since the new method greatly enlarges the convergence domain, the number of initial values of hypocenter that we need to select in solving the earthquake location problem can be significantly reduced compared to the conventional method. This greatly reduces the overall computational cost.

\section{Numerical Experiments} \label{sec:num}
In this section, three examples are presented to demonstrate the validity of our method. And we will see the comparison between the conventional method and the new method for the earthquake location problem.

\begin{figure*} 
	\begin{tabular}{cc}
		\includegraphics[width=0.38\textwidth, height=0.15\textheight]{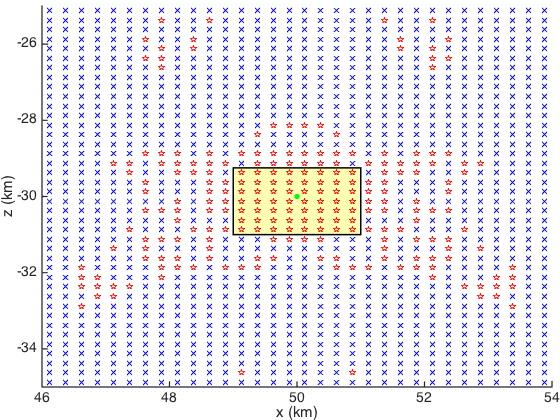} &
		\includegraphics[width=0.38\textwidth, height=0.15\textheight]{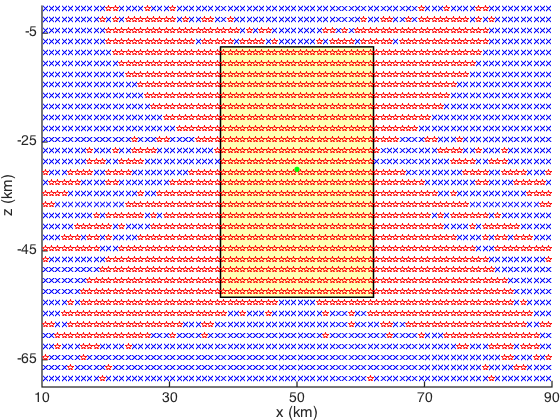} \\
		\includegraphics[width=0.38\textwidth, height=0.15\textheight]{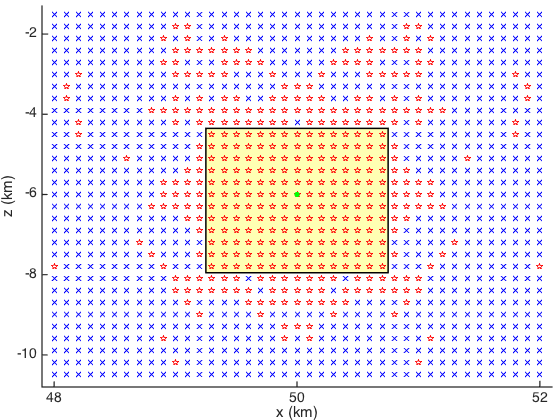} &
		\includegraphics[width=0.38\textwidth, height=0.15\textheight]{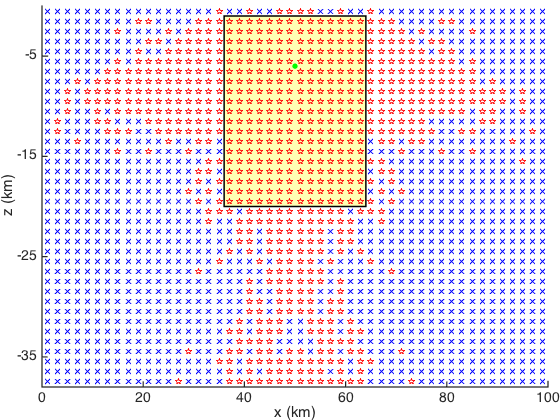}
	\end{tabular}
	\caption{Illustration of the Exam \ref{exam:homo}. The green point is the real hypocenter. The red pentagram and the blue x-mark indicate the initial hypocenter at this location converge and misconvergence to the real hypocenter respectively. Left: the conventional method; Right: the new method. Up figures for deep earthquake $\bxi_T=(50km,-30km)$ and Bottom figures for shallow earthquake $\bxi_T=(50km,-6km)$. In the light yellow rectangular region, all the tested initial hypocenter converge to the correct hypocenter.} \label{fig:exam41}
\end{figure*}

\begin{example} \label{exam:homo}
	Let's take the same parameters set up as in Exam \ref{exam:large_pert}. Then the real signals $d_r(t)$ and the synthetic $s(\boeta_r,t)$ can be obtained by \eqref{eqn:d_rel}, \eqref{eqn:s_rel} and \eqref{sol:analy_2D} for different receiver $r=1,2,\cdots,20$. 
	
	Consider an earthquake occurs at hypocenter $\bxi_T=(50km,-30km)$ and origin time $\tau_T=10s$ with dominant frequency $f_0=2$Hz. In Figure \ref{fig:exam41} (Up, Left), $1280$ uniformly distributed grid nodes are tested as the initial hypocenter of earthquake $\bxi$ in the searching domain $[46km,54km]\times[-35km,-25km]$ for the conventional method. This is the simplest experimental design technique \cite{CoCo:92}, but it is enough to illustrate the effectiveness of our new method. There are $228$ grid nodes converge to the correct hypocenter. Therefore, the area of the convergence domain is roughly estimated as
	\begin{equation*}
		(54-46)\times((-25)-(-35))\times\frac{228}{1280}=14.3 \;km^2.
	\end{equation*}
	
	In Figure \ref{fig:exam41}  (Up, Right), $2800$ uniformly distributed grid nodes are tested as the initial hypocenter of earthquake $\bxi$ in the searching domain $[10km,90km]\times[-70km,0km]$ for the new method. There are $1597$ grid nodes converge to the correct hypocenter. Using the same formula, the area of the convergence domain is roughly estimated as $3194\;km^2$. In contrast, the convergence probability of the new method is about
	\begin{equation*}
		3194\div 14.3\approx 223,
	\end{equation*}
	times that of the conventional method for this case.
	
	From the figure, we can also see that all the tested initial hypocenter in the rectangular region $[49km,51km]\times[-31km,-29.25km]$ converge to the correct hypocenter for the conventional method. For the new method, this rectangular region is $[38km,62km]\times[-53.5km,-7.5km]$, its area is $315$ times of the former for this case.

	Consider an earthquake occurs at hypocenter $\bxi_T=(50km,-6km)$ and origin time $\tau_T=10s$ with dominant frequency $f_0=2$Hz. In Figure \ref{fig:exam41} (Bottom, Left), $1271$ uniformly distributed grid nodes are tested as the initial hypocenter of earthquake $\bxi$ in the searching domain $[48km,52km]\times[-11km,0km]$ for the conventional method. There are $411$ grid nodes converge to the correct hypocenter, thus the area of the convergence domain is roughly estimated as $14.2\;km^2$. In Figure \ref{fig:exam41}  (Bottom, Right), $1900$ uniformly distributed grid nodes are tested as the initial hypocenter of earthquake $\bxi$ in the searching domain $[0km,100km]\times[-38km,0km]$ for the new method. There are $740$ grid nodes converge to the correct hypocenter, thus the area of the convergence domain is roughly estimated as $1480\;km^2$. In contrast, the convergence probability of the new method is about $104$ times that of the conventional method for this case.

	From the figure, we can also see that all the tested initial hypocenter in the rectangular region $[49.25km,50.75km]\times[-7.95km,-4.35km]$ converge to the correct hypocenter for the conventional method. For the new method, this rectangular region is $[36km,64km]\times[-20km,-1km]$, its area is $98$ times of the former for this case. 
	
	Considering all of the above, we note that the new method works better for the deep earthquake rather than the shallow earthquake. One explanation is that the convergence domain is nearly symmetric about the earthquake hypocenter. But it doesn't hold for shallow earthquake in $z$ direction since selecting the initial hypocenter above the surface is non-physical. This discussion also applies to the following examples. $\Box$
\end{example}

\begin{figure} 
	\centering
	\includegraphics[width=0.60\textwidth, height=0.25\textheight]{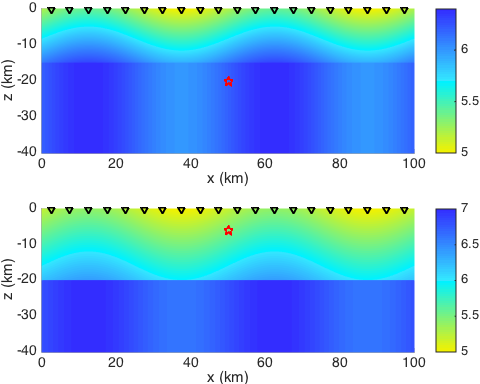}
	\caption{Velocity models in Example \ref{exam:two_layer}. The red pentagrams show the hypocenter of earthquake and the black triangles indicate the receiver stations.} \label{fig:velocity_illu}
\end{figure}

\begin{figure*} 
	\begin{tabular}{cc}
		\includegraphics[width=0.38\textwidth, height=0.15\textheight]{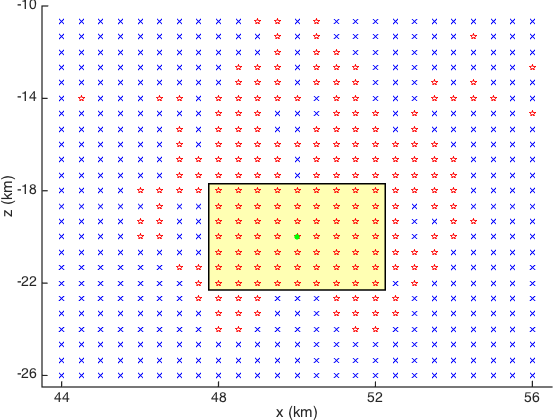} &
		\includegraphics[width=0.38\textwidth, height=0.15\textheight]{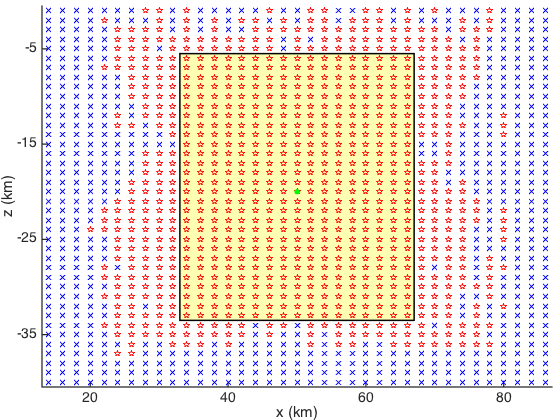} \\
		\includegraphics[width=0.38\textwidth, height=0.15\textheight]{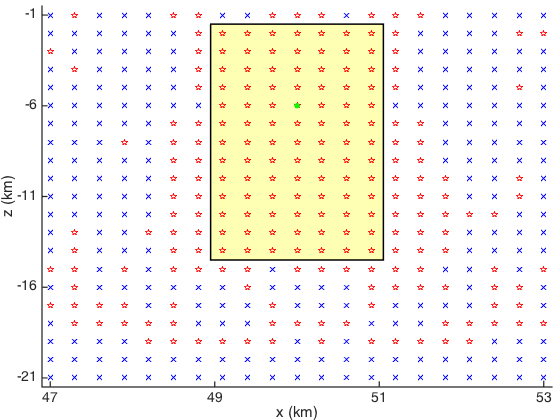} &
		\includegraphics[width=0.38\textwidth, height=0.15\textheight]{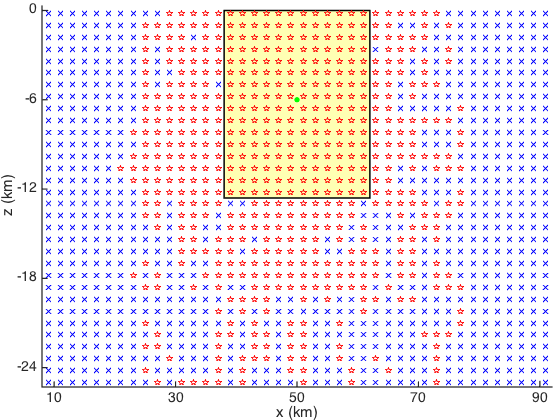}
	\end{tabular}
	\caption{Illustration of the Exam \ref{exam:two_layer}. The green point is the real hypocenter. The red pentagram and the blue x-mark indicate the initial hypocenter at this location converge and misconvergence to the real hypocenter respectively. Left: the conventional method; Right: the new method. Up figures for deep earthquake $\bxi_T=(50km,-20km)$ and Bottom figures for shallow earthquake $\bxi_T=(50km,-6km)$. In the light yellow rectangular region, all the tested initial hypocenter converge to the correct hypocenter.} \label{fig:exam42}
\end{figure*}

\begin{example} \label{exam:two_layer}
	Consider the two-layer model in the bounded domain $[0km,100km]\times[-40km,0km]$, the wave speed is 
	\begin{equation*}
		c(x,z)=\left\{\begin{array}{ll}
			5.2-0.06z+0.2\sin\frac{\pi x}{25}, & -15km\le z\le 0km, \\
			6.2+0.2\sin\frac{\pi x}{25}, & -40km\le z<-15km,
		\end{array}\right.
	\end{equation*}	
	for depth earthquake and
	\begin{equation*}
		c(x,z)=\left\{\begin{array}{ll}
			5.2-0.05z+0.2\sin\frac{\pi x}{25}, & -20km\le z\le 0km, \\
			6.8+0.2\sin\frac{\pi x}{25}, & -40km\le z<-20km,
		\end{array}\right.
	\end{equation*}		
	for shallow earthquake. The unit is `km/s'. We use the finite difference scheme \cite{Da:86,WaYaSo:12,YaTeZhLi:03,YaLuWuPe:04,YaPeLuTe:06} to solve the acoustic wave equation \eqref{eqn:wave} with initial condition \eqref{ic:wave}. The problem can also be solved by other numerical methods, e.g. finite element methods \cite{LyDr:72,Ma:84}, the spectral element method \cite{KiLiTr:11,LiPoKoTr:04} and the discontinuous Galerkin method \cite{HeYaWu:14,HeYaWu:15}. The reflection boundary condition is used on the earth's surface, and the perfectly matched layer \cite{KoTr:03,MaYaSo:15} is used for other boundaries. The delta function $\delta(\bx-\bxi)$ in the wave equation \eqref{eqn:wave} is discretized using the techniques proposed in \cite{We:08}.
	\begin{equation*}
		\delta(x)=\left\{\begin{array}{ll}
			\frac{1}{h}\left(1-\frac{5}{4}\abs{\frac{x}{h}}^2-\frac{35}{12}\abs{\frac{x}{h}}^3
				+\frac{21}{4}\abs{\frac{x}{h}}^4-\frac{25}{12}\abs{\frac{x}{h}}^5\right), & \abs{x}\le h, \\
			\frac{1}{h}\left(-4+\frac{75}{4}\abs{\frac{x}{h}}-\frac{245}{8}\abs{\frac{x}{h}}^2+\frac{545}{24}\abs{\frac{x}{h}}^3
				-\frac{63}{8}\abs{\frac{x}{h}}^4+\frac{25}{24}\abs{\frac{x}{h}}^5\right), & h<\abs{x}\le 2h, \\
			\frac{1}{h}\left(18-\frac{153}{4}\abs{\frac{x}{h}}+\frac{255}{8}\abs{\frac{x}{h}}^2-\frac{313}{24}\abs{\frac{x}{h}}^3
				+\frac{21}{8}\abs{\frac{x}{h}}^4-\frac{5}{24}\abs{\frac{x}{h}}^5\right), & 2h<\abs{x}\le 3h, \\
			0, & \abs{x}>3h.
		\end{array}\right. 
	\end{equation*}
	There are 20 equidistant receivers on the surface
	\begin{equation*}
		\boeta_r=(x_r,z_r)=(5r-2.5km,0), \quad r=1,2,\cdots,20.
	\end{equation*}
	Since the hypocenter of earthquake is not far from the receiver stations, we only use the direct wave to locate the earthquake.
	
	Consider an earthquake occurs below the medium interface $\bxi_T=(50km,-20km)$ and origin time $\tau_T=10s$ with dominant frequency $f_0=2$Hz (see Figure \ref{fig:velocity_illu} Up). In Figure \ref{fig:exam42} (Up, Left), $600$ uniformly distributed grid nodes are tested as the initial hypocenter of earthquake $\bxi$ in the searching domain $[44km,56km]\times[-26km,-10km]$ for the conventional method. There are $210$ grid nodes converge to the correct hypocenter, thus the area of the convergence domain is roughly estimated as $67.2\;km^2$. In Figure \ref{fig:exam42}  (Up, Right), $1480$ uniformly distributed grid nodes are tested as the initial hypocenter of earthquake $\bxi$ in the searching domain $[10km,90km]\times[-40km,0km]$ for the new method. There are $881$ grid nodes converge to the correct hypocenter, thus the area of the convergence domain is roughly estimated as $1905\;km^2$. In contrast, the convergence probability of the new method is about $28.3$ times that of the conventional method for this case.
	
	From the figure, we can also see that all the tested initial hypocenter in the rectangular region $[48.95km,51.05km]\times[-14.5km,-1.5km]$ converge to the correct hypocenter for the conventional method. For the new method, this rectangular region is $[33km,67km]\times[-33.5km,-5.5km]$, its area is $46.0$ times of the former for this case.
	
	Consider an earthquake occurs above the medium interface $\bxi_T=(50km,-6km)$ and origin time $\tau_T=10s$ with dominant frequency $f_0=2$Hz (see Figure \ref{fig:velocity_illu} Bottom). In Figure \ref{fig:exam42} (Bottom, Left), $441$ uniformly distributed grid nodes are tested as the initial hypocenter of earthquake $\bxi$ in the searching domain $[47km,53km]\times[-21km,-1km]$ for the conventional method. There are $216$ grid nodes converge to the correct hypocenter, thus the area of the convergence domain is roughly estimated as $58.8\;km^2$. In Figure \ref{fig:exam42}  (Bottom, Right), $1344$ uniformly distributed grid nodes are tested as the initial hypocenter of earthquake $\bxi$ in the searching domain $[8km,92km]\times[-25km,0km]$ for the new method. There are $592$ grid nodes converge to the correct hypocenter, thus the area of the convergence domain is roughly estimated as $925\;km^2$. In contrast, the convergence probability of the new method is about $15.7$ times that of the conventional method for this case.
	
	From the figure, we can also see that all the tested initial hypocenter in the rectangular region $[47.75km,52.25km]\times[-22.3km,-17.7km]$ converge to the correct hypocenter for the conventional method. For the new method, this rectangular region is $[38km,62km]\times[-12.6km,0km]$, its area is $11.1$ times of the former for this case. $\Box$
\end{example}

\begin{figure} 
	\centering
	\includegraphics[width=0.60\textwidth, height=0.25\textheight]{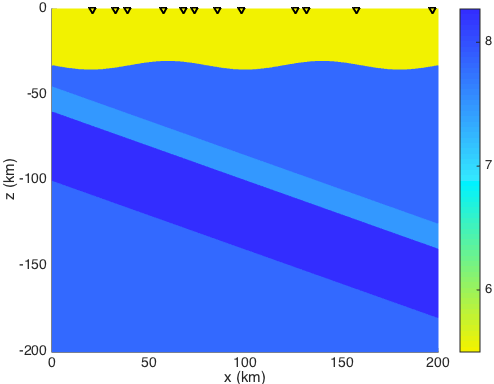}
	\caption{Velocity model in Example \ref{exam:practical}. The black triangles indicate the receiver stations.} \label{fig:pracvel_illu}
\end{figure}

\begin{figure*} 
	\centering
	\includegraphics[width=0.76\textwidth, height=0.30\textheight]{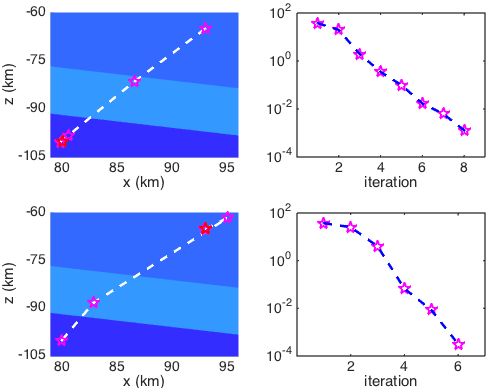}
	\caption{Convergent history of the first case in Example \ref{exam:practical}, from initial hypocenter in the subduction zone to the real hypocenter in the mantle (Up) and its contrary case (Bottom). Left: the convergent trajectories; Right: the absolute errors with respect to iteration step between the real and computed hypocenter of the earthquake.} \label{fig:exam43_man_sub}
\end{figure*}

\begin{figure*} 
	\centering
	\includegraphics[width=0.76\textwidth, height=0.16\textheight]{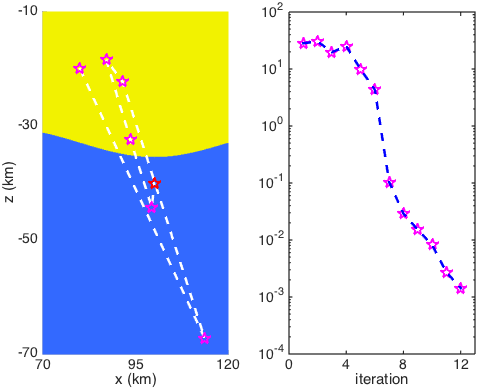}
	\caption{Convergent history of the second case in Example \ref{exam:practical},  from initial hypocenter in the crust to the real hypocenter in the mantle. Left: the convergent trajectory; Right: the absolute errors with respect to iteration step between the real and computed hypocenter of the earthquake.} \label{fig:exam43_man_crust}
\end{figure*}

\begin{example} \label{exam:practical}
	Consider the velocity model consisting of the crust and the mantle, containing an undulated Moho discontinuity and a subduction zone with a thin low velocity layer atop a fast velocity layer \cite{ToYaLiYaHa:16}, see Figure \ref{fig:pracvel_illu} for illustration. The computational domain is $[0km,200km]\times[-200km,0km]$, and the wave speed is 
	\begin{equation*}
		c(x,z)=\left\{\begin{array}{ll}
			5.5, & -33-2.5\sin\frac{\pi x}{40}\le z<0, \\
			7.8, & -45-0.4x\le z<-33-2.5\sin\frac{\pi x}{40}, \\
			7.488, & -60-0.4x\le z<-45-0.4x. \\
			8.268, & -100-0.4x\le z<-60-0.4x. \\
			7.8, & \textnormal{others}.
		\end{array}\right.
	\end{equation*}	
	with unit `km/s'. We consider the same set-up as in Exam \ref{exam:two_layer}, e.g. the forward scheme, the boundary conditions and the discretized delta function. There are 12 receivers $\boeta_r=(x_r,z_r)$ on the surface with $z_r=0$, their horizontal positions are randomly given, see Table \ref{tab:pracvel_recepos} for details. In this example, we still only use the directive wave to locate the earthquake. In real world, region with the similar velocity model is always seismogenic zone \cite{ToZhYa:11}. Earthquakes in this kind of region can occur in the crust, in the subduction zone or in the mantle \cite{ToZhYa:12}. Complex velocity structure makes source location very difficult.
	
	We firstly investigate the case that the earthquake occurs in the mantle but the initial hypocenter of the earthquake is chosen in the subduction zone, and its contrary case. In Figure \ref{fig:exam43_man_sub}, we can see the convergent history. The second case is that the earthquake occurs in the mantle but the initial hypocenter of the earthquake is chosen in the crust. The convergent history can be seen in Figure \ref{fig:exam43_man_crust}. From these tests, we can observes nice convergent result of the new method, even though the real and initial hypocenter of the earthquakes are far from each other. $\Box$
\end{example}

\begin{table*}
    	\caption{Example \ref{exam:practical}: the horizontal positions of receivers, with unit `km'.} \label{tab:pracvel_recepos}
	\begin{center}\begin{tabular}{ccccccccccccc} \hline
		$r$ & $1$ & $2$ & $3$ & $4$ & $5$ & $6$ & $7$ & $8$ & $9$ & $10$ & $11$ & $12$ \\ \hline
		$z_r$ & $21$ & $33$ & $39$ & $58$ & $68$ & $74$ & $86$ & $98$ & $126$ & $132$ & $158$ & $197$ \\ \hline
    	\end{tabular}\end{center}
\end{table*}

\section{Conclusion and Discussion} \label{sec:con}
The main contribution in this paper is that convergence domain of the waveform based earthquake location method has been greatly expanded. Accordingly to the numerical evidence presented earlier, the convergence domain has been enlarged $10\sim 300$ times in the two test problems. This means that even from the relatively poor initial values of earthquake hypocenter, our method is also likely to convergence to the correct results with high accuracy.

We have to explain that this paper focuses on the development of new method. For practical three-dimensional problem, we believe that the new method is also applicable. We are investigating this approach. We hope this can be reported in an independent publication in the near future.

\begin{acknowledgments}
This work was supported by the National Nature Science Foundation of China (Grant Nos 41230210, 41390452). Hao Wu was also partially supported by the National Nature Science Foundation of China (Grant No 11101236) and SRF for ROCS, SEM. The authors are grateful to Prof. Shi Jin for his helpful suggestions and discussions that greatly improve the presentation. Hao Wu would like to thank Prof. Ping Tong for his valuable comments.
\end{acknowledgments}

\label{lastpage}

\end{document}